\documentclass[12pt]{amsart}

\usepackage{amsmath,amssymb}
\newtheorem{thm}{Theorem}
\newtheorem{lem}[thm]{Lemma}
\newtheorem{prop}[thm]{Proposition}
\newtheorem{cor}[thm]{Corollary}

\newcommand{\fq}{\mathbb{F}_q}

\numberwithin{equation}{section}

\newcommand{\Z}{{\mathbb Z}} 

\newcommand{\R}{{\mathbb R}}

\newcommand{\FF}{{\mathbb F}}

\renewcommand{\mod}{\;\operatorname{mod}}

\newcommand{\tr}{\operatorname{tr}}
\newcommand{\intinf}{\int_{-\infty}^\infty}
\newcommand{\ave}[1]{\left\langle#1\right\rangle} 

\newcommand{\leg}[2]{\left( \frac{#1}{#2} \right)}  

\newcommand{\USp}{\operatorname{USp}}

\newcommand{\LL}{\mathcal L} 

\begin{document}
\baselineskip=17pt  

\title[Traces of high powers]
{Traces of high powers of the Frobenius class in the hyperelliptic ensemble} 

\author{Ze\'ev Rudnick}
\address{Raymond and Beverly Sackler School of Mathematical Sciences,
Tel Aviv University, Tel Aviv 69978, Israel and 
School of Mathematics, Institute for Advanced Study,   
Princeton, NJ 08540 }
\email{rudnick@post.tau.ac.il}

\date{September 1, 2009}




 \begin{abstract}
 
 The zeta function of a curve over a finite field may be expressed in terms of 
the characteristic polynomial of a unitary symplectic matrix $\Theta$, 
called the Frobenius class of the curve. We compute the expected value of $\tr (\Theta^n)$ 
for an ensemble of hyperelliptic curves of genus $g$ over a fixed finite field in the limit of large genus, 
and compare the results to the corresponding averages over the unitary symplectic group $\USp(2g)$. 
We are able to compute the averages for powers $n$ almost up to $4g$, finding agreement with the Random Matrix 
results except for small $n$ and for $n=2g$. 
As an application we compute the one-level density of zeros of the zeta function of the curves,  
including lower-order terms, for test  functions whose Fourier transform is supported in $(-2,2)$. 
The results confirm in part a conjecture of Katz and Sarnak, that to leading order the low-lying zeros 
for this ensemble have symplectic statistics. 
 
 \end{abstract}

\maketitle

\section{Introduction} 
 
Fix a finite field $\fq$ of odd cardinality, and let $C$ be 
a non singular projective   curve defined over $\fq$.  
For each extension field of degree $n$ of $\fq$, denote by $N_n(C)$ 
the number of points of $C$ in $\mathbb F_{q^n}$. 
The zeta function associated to $C$ is defined as 
$$Z_C(u) = \exp \sum_{n=1}^\infty N_n(C)\frac{u^n}n, \quad |u|<1/q$$
and is known to be a rational function of $u$ of the form
\begin{equation}\label{zeta function}
 Z_C(u) = \frac{P_C(u)}{(1-u)(1-qu)} 
\end{equation}
where $P_C(u)$ 
is a polynomial of degree $2g$ with integer coefficients, 
satisfying a functional equation
$$ P_C(u) = (q u^2)^g P_C(\frac 1{qu})\;.
$$
The Riemann Hypothesis, proved by Weil \cite{Weil}, is that the
zeros of $P(u)$ all lie on the circle $|u|=1/\sqrt{q}$. 
Thus one may give a spectral interpretation of $P_C(u)$ as the characteristic
polynomial of a $2g\times 2g$ unitary 
matrix $\Theta_C$
$$ P_C(u) = \det (I-u\sqrt{q}\Theta_C)$$
so that the eigenvalues $e^{i\theta_j}$ of $\Theta_C$ correspond to
zeros $q^{-1/2}e^{-i\theta_j}$ of $Z_C(u)$. 
The matrix (or rather the conjugacy class) $\Theta_C$ is called the
unitarized Frobenius class of $C$. 

We would like to study the how the Frobenius classes $\Theta_C$ change
as we vary the curve over a family of hyperelliptic curves of genus
$g$, in the limit of large genus and fixed constant field. 
The particular family $\mathcal H_{2g+1}$   we choose is the family
of all curves given in affine form by an equation
$$ C_Q: y^2 = Q(x)$$
where
$$ Q(x) = x^{2g+1} +a_{2g} +\dots +a_0 \in \fq[x]$$
is a squarefree, monic polynomial of degree $2g+1$.
The curve $C_Q$ is thus nonsingular and of genus $g$.

We consider $\mathcal H_{2g+1}$ as a probability space (ensemble)
with the uniform probability measure, so that the expected value of
any function $F$ on $\mathcal H_{2g+1}$ is defined as 
\begin{equation*}
\ave{F}:=\frac 1{\#\mathcal H_{2g+1}} \sum_{Q\in \mathcal H_{2g+1}} F(Q)
\end{equation*}

Katz and Sarnak showed \cite{KS} that as $q\to \infty$, the Frobenius classes $\Theta_Q$ become
equidistributed in the unitary symplectic group $\USp(2g)$ 
(in genus one this is due to Birch \cite{Birch} for $q$ prime,   
and to Deligne \cite{Deligne}).  
That is for any continuous function on 
the space of conjugacy classes of $\USp(2g)$, 
$$\lim_{q\to\infty}  \ave{F(\Theta_Q)} = \int_{\USp(2g)} F(U)dU $$
This implies that various statistics of the eigenvalues can, in this limit,
be computed by integrating the corresponding quantities over $\USp(2g)$.

Our goal is to explore the opposite limit, that of fixed constant
field and large genus ($q$ fixed, $g\to \infty$), cf \cite{KR08, FR}. Since the matrices 
$\Theta_Q$ now inhabit different  spaces as $g$ grows,
it is not clear how to formulate an equidistribution problem.
However one can still meaningfully discuss various statistics, 
the most fundamental being various products  of traces of powers of
$\Theta_Q$, that is $\ave{\prod_{j=1}^r \tr (\Theta_Q^{n_j}) }$. 
Here we study the basic case of the expected values 
$\ave{\tr  \Theta_Q^n}$ where $n$ is of order of the genus $g$.

The mean value of traces of powers when averaged over the 
unitary symplectic group $\USp(2g)$ are easily seen to be  \cite{DS} 
\begin{equation}\label{traces in RMT}
\int_{\USp(2g)} \tr (U^n) dU  =
\begin{cases}
2g& n=0 \\ 
-\eta_n & 1\leq |n|\leq 2g \\ 
0& |n|>2g          
\end{cases}
\end{equation}
where
$$ \eta_n= \begin{cases}
           1& n \mbox{ even}\\0& n\mbox{ odd} 
           \end{cases} 
$$

We will show: 
\begin{thm}\label{one level} 
For all $n>0$ we have 
\begin{equation*}
\begin{split}
\ave{\tr \Theta_Q^n} &= 
\left\{\begin{matrix} -\eta_n, & 0<n<2g\\  \\ -1-\frac 1{q-1},&
n=2g\\  \\0,& n>2g\end{matrix} \right \}  
+\eta_n\frac 1{q^{n/2}} \sum_{\substack{\deg P \mid \frac n2\\P \mbox{ prime}}} \frac{\deg P}{|P|+1}
\\
&+O_q\left(n q^{n/2-2g} +g q^{-g} \right)
\end{split}
\end{equation*}
the sum over all irreducible monic polynomials $P$, and
where $|P|:=q^{\deg P}$. 
\end{thm}
In particular we have
\begin{cor}
If   $3\log_qg< n< 4g-5\log_q g$  but $n\neq 2g$ 
then 
$$ 
\ave{\tr \Theta_Q^n}  =\int_{\USp(2g)} \tr U^n dU +o(\frac 1g) \;.
$$
\end{cor}
We do however get deviations from the Random Matrix Theory results \eqref{traces in RMT} 
for small values of $n$, for instance 
 $$ \ave{\tr \Theta_Q^2} \sim \int_{\USp(2g)} \tr U^{2} dU  + \frac 1{q+1} $$
and for $n=2g$ where we have 
$$ 
\ave{\tr \Theta_Q^{2g}}  \sim \int_{\USp(2g)} \tr U^{2g} dU -\frac 1{q-1} \;.
$$



Analogous results can be derived for mean values of products, e.g. for  
$\ave{\tr \Theta_Q^m\tr  \Theta_Q^n}$, when $m+n<4g$, see \S~\ref{sec:products}.

To prove these results, we cannot use the powerful equidistribution theorem of Deligne \cite{Deligne}, 
as was done for the fixed genus case in \cite{KS}. Rather, we use a variant of the analytic methods 
developed to deal with such problems in the number field setting \cite{OS, Sarnakappendix, Sound}.   
Extending the range of our results to cover $n>4g$ is a challenge.

\subsection{Application: The one-level density} 

The traces of powers determine all {\em linear} statistics,    
such as the number of angles $\theta_j$ lying in
a subinterval of $\R/2\pi \Z$, or the one-level density,
a smooth linear statistic. To define the one-level density, 
we start with an even test function $f$, say in the Schwartz space $\mathcal S(\R)$, 
and for any $N\geq 1$ set 
\begin{equation*}
 F(\theta):=\sum_{k\in \Z} f(N(\frac {\theta}{2\pi} -k))
\end{equation*}
which has period $2\pi$ and is localized in an interval of size
$\approx 1/N$ in $\R/2\pi\Z$.  
For a unitary $N\times N$ matrix $U$ with eigenvalues $e^{ i \theta_j} $, $j=1,\dots N$, define 
$$ Z_{f}(U):= \sum_{j=1}^N F(\theta_j) $$
which counts the number of ``low-lying'' eigenphases $\theta_j$ in the smooth
interval of length $\approx 1/N$ around the origin defined by $f$.

Katz and Sarnak conjectured \cite{KS-BAMS} that for fixed $q$, the
expected value of $Z_f$ over $\mathcal H_{2g+1}$  
will converge to $\int_{\USp(2g)}Z_f(U)dU $ as $g\to \infty$ 
for {\em any} such test function $f$. 
Theorem~\ref{one level}  implies: 
\begin{cor}\label{cor Zf}
If $f\in \mathcal S(\R)$ is even, with Fourier transform $\^f$ 
supported in $(-2,2)$ then   
$$ \ave{Z_{f}} = \int_{\USp(2g)} Z_f(U)dU + \frac{\mbox{dev}(f)}{g}
+o(\frac 1g) 
$$
where 
$$dev(f) = \^f(0)\sum_{P \mbox{ prime} } \frac{\deg P}{|P|^2-1 }
-\^f(1) \frac 1{q-1} 
$$
the sum over all irreducible monic polynomials $P$. 
\end{cor}
To show corollary~\ref{cor Zf}, one uses a Fourier expansion 
to see that 
\begin{equation}\label{fourier for nu}
Z_f(U) = 
\intinf f(x)dx +\frac1N\sum_{n\neq 0} \^f(\frac nN) \tr U^n \;.
\end{equation}
Averaging $Z_f(U)$ over the symplectic group $\USp(2g)$, 
using \eqref{traces in RMT}, and assuming $f$ is even, 
gives\footnote{Note that as $g\to \infty$, 
$\int_{\USp(2g)}Z_f(U)dU \sim 
\intinf f(x) \left(1-\frac{\sin 2\pi x}{2\pi x} \right) dx $}
\begin{equation*}
\int_{\USp(2g)}Z_f(U)dU =\^f(0) -\frac 1{g}\sum_{1\leq m\leq g} \^f(\frac mg) 
\end{equation*}
and then we use Theorem~\ref{one level} to deduce Corollary~\ref{cor Zf}.

Corollary~\ref{cor Zf} is completely analogous to what is known in the number field setting for the corresponding case of 
zeta functions of quadratic fields,   
except for the lower order term which is different: 
While the coefficient of $\^f(0)$ is as in the number field setting \cite{SJMiller lower}, 
the coefficient of $\^f(1)$ is special to our function-field setting.  
 
\subsection{Acknowledgments} 
We thank P\"ar Kurlberg and Nicolas Templier for their comments. 
Supported by the Israel Science Foundation (grant No. 925/06) and by the Oswald Veblen Fund at the Institute for Advanced Study, Princeton.

\section{Quadratic L-functions} 
In this section we give some known background on the zeta function of 
hyperelliptic curves. 
The theory was initiated by E. Artin \cite{Artin}. 
We use Rosen \cite{Rosen} as a general reference. 

\subsection{} 
For a nonzero polynomial $f\in \fq[x]$, we define the norm $|f|:=q^{\deg f}$. 
A ``prime'' polynomial is a monic irreducible polynomial. 
For a monic polynomial $f$, The von Mangoldt function $\Lambda(f)$ is
defined to be zero unless $f=P^k$  is a prime power in which case
$\Lambda(P^k)=\deg P$.  

The analogue of Riemann's zeta function is 
$$\zeta_q(s):= \prod_{P \mbox{ prime}} (1-|P|^{-s})^{-1}$$
which is shown to equal 
\begin{equation}\label{simple zeta}
  \zeta_q(s) = \frac 1{1-q^{1-s}}
\end{equation}

Let $\pi_q(n)$ be the number of prime polynomials of degree $n$. 
The Prime Polynomial Theorem in $\fq[x]$  asserts that  
$$ \pi_q(n) = \frac{q^n}n +O(q^{n/2})$$
which follows from the identity (equivalent to \eqref{simple zeta})
\begin{equation}\label{equiv simple zeta}
  \sum_{\deg(f)=n} \Lambda(f) = q^n 
\end{equation}
the sum over all monic polynomials of degree $n$.

\subsection{} 
For a monic polynomial  $D\in \FF_q[x]$ of positive degree, which is not a perfect square,  
we define the quadratic character $\chi_D$ in terms of the quadratic residue 
symbol for $\fq[x]$  by 
$$ \chi_D(f) = \leg{D}{f}$$
and the corresponding L-function
$$
\LL(u,\chi_D):= \prod_P (1-\chi_D(P)u^{\deg P})^{-1}, \quad |u|<\frac 1q 
$$
the product over all monic irreducible (prime) polynomials $P$. 
Expanding in additive form using unique factorization, we write 
$$\LL(u,\chi_D) = \sum_{\beta\geq 0} A_D(\beta)u^\beta
$$
with 
$$A_D(\beta):= \sum_{\substack{\deg B=\beta\\B \mbox{ monic}}} \chi_D(B) \;.
$$
If $D$ is non-square  of positive degree, then $A_D(\beta) =  0$ 
for $\beta \geq \deg D$
and hence the L-function is in fact a polynomial of degree at most
$\deg D-1$. 

\subsection{} 
To proceed further,  assume that $D$ is square-free 
(and monic of positive degree).   
Then $\LL(u,\chi_D)$ has a ``trivial'' zero at $u=1$  if and only if $\deg D$ is even. Thus  
$$\LL(u,\chi_D)=(1-u)^\lambda \LL^*(u,\chi_D), \quad 
\lambda=
\begin{cases} 1& \deg D \mbox{ even}\\0& \deg(D)\mbox{ odd}\end{cases}
$$
where $\LL^*(u,\chi_D)$ is a polynomial of even degree 
$$2\delta=\deg D-1-\lambda$$
satisfying the functional equation  
\begin{equation*}
  \LL^*(u,\chi_D) = (qu^2)^{\delta} \LL^*(\frac 1{qu},\chi_D) \;.
\end{equation*}
In fact $\LL^*(u,\chi_D)$ is the Artin L-function associated to the unique nontrivial quadratic character of $\fq(x)(\sqrt{D(x)})$, 
see \cite[Propositions 17.7 and 14.6]{Rosen}. 
We write  

$$ \LL^*(u,\chi_D) = \sum_{\beta=0}^{2\delta} A^*_D(\beta) u^\beta$$
where $A_D^*(0)=1$, and the coefficients $A_D^*(\beta)$  satisfy
\begin{equation}\label{FE for coeffs}
 A_D^*(\beta) = q^{\beta-\delta}A_D^*(2\delta-\beta) \;.
\end{equation}
In particular the leading coefficient is  $A_D^*(2\delta) = q^\delta$.

\subsection{} 
For $D$ monic, square-free, and of positive degree, 
the zeta function \eqref{zeta function} of the hyperelliptic curve $y^2=D(x)$
is 
$$Z_D(u) = \frac{\LL^*(u,\chi_D)}{(1-u)(1-qu)}  \;.$$
The Riemann Hypothesis, proved by Weil \cite{Weil}, asserts that all zeros of $Z_C(u)$, 
hence of $\LL^*(u,\chi_D)$, lie on the circle $|u|=1/\sqrt{q}$. 
Thus we may write 
$$\LL^*(u,\chi_D) = \det(I-u\sqrt{q} \Theta_D) $$
for a unitary $2\delta\times 2\delta$ matrix $\Theta_D$.

\subsection{} 

By taking a logarithmic derivative of the identity 
$$
\det(I-u\sqrt{q}\Theta_D) = (1-u)^{-\lambda} \prod_P (1-\chi_D(P)u^{\deg P})^{-1}
$$
which comes from writing $\LL^*(u,\chi_D) = (1-u)^{-\lambda} \LL(u,\chi_D)$, 
we find 
\begin{equation}\label{explicit formula} 
-  \tr \Theta_D^n =\frac  \lambda{q^{n/2}} + 
\frac 1{q^{n/2}} \sum_{\deg f=n} \Lambda(f) \chi_D(f)
\end{equation}

\subsection{}
Assume now that $B$ is monic, of positive degree and not a perfect square. 
Then we have a bound for the character sum over primes: 
\begin{equation}\label{estimate on prime sum} 
  \left| \sum_{\substack{\deg P=n\\P \mbox{ prime}}} 
 \leg{B}{P}\right| \ll \frac{\deg B}{n} q^{n/2}
\end{equation}
This is deduced by writing $B=DC^2$ with $D$  square-free, 
of positive degree, and then using the explicit formula 
\eqref{explicit formula}  
and the unitarity of $\Theta_D$  (which is the Riemann Hypothesis). 

\section{The hyperelliptic ensemble $\mathcal H_{2g+1}$} 
\subsection{Averaging over $\mathcal H_{2g+1}$} 

We denote by $\mathcal H_d$ the set of square-free monic polynomials
of degree $d$ in $\fq[x]$. 
The cardinality of $\mathcal H_{d}$ is
\begin{equation*}
\# \mathcal H_{d} = \begin{cases}
(1-\frac 1q)q^{d},& d\geq 2\\q,&d=1
\end{cases}
\end{equation*} 
as is seen by writing 
$$\sum_{d\geq 0}\frac{ \#\mathcal H_{d}}{ q^{ds}} =\sum_{f \mbox{ monic
    squarefree}} |f|^{-s} = \frac{\zeta_q(s)}{\zeta_q(2s)}$$ 
and using \eqref{simple zeta}.  
In particular for $g\geq 1$, 
$$\#\mathcal H_{2g+1} = (q-1)q^{2g} \;.
$$

We consider $\mathcal H_{2g+1}$ as a probability space (ensemble)
with the uniform probability measure, so that the expected value of
any function $F$ on $\mathcal H_{2g+1}$ is defined as 
\begin{equation}
\ave{F}:=\frac 1{\#\mathcal H_{2g+1}} \sum_{Q\in \mathcal H_{2g+1}} F(Q)
\end{equation}

We can pick out square-free polynomials by using the
M\"obius function $\mu$ of $\fq[x]$ (as is done over the integers) via 
$$ \sum_{A^2\mid Q} \mu(A) = 
\begin{cases}
  1& Q \mbox{ square-free}\\0& \mbox{ otherwise}
\end{cases}
$$
Thus we may write  expected values as  
\begin{equation}\label{ave with mobius}
\ave{F(Q)} = \frac 1 {(q-1)q^{2g}} \sum_{2\alpha+\beta=2g+1} 
\sum_{\deg B=\beta}\sum_{\deg A=\alpha} \mu(A) F(A^2B)
\end{equation}
the sum over all monic $A$, $B$.

\subsection{Averaging quadratic characters} 
Suppose now that we are given a polynomial $f\in \fq[x]$ and apply \eqref{ave with mobius} 
to the quadratic character $\chi_Q(f)=\leg{Q}{f}$. 
Then 
$$\chi_{A^2B}(f) = \leg{B}{f} \leg{A}{f}^2 = 
\begin{cases}
  \leg{B}{f}& \gcd(A,f)=1 \\0& \mbox{ otherwise}
\end{cases}
$$
Hence 
\begin{equation*}
  \ave{\chi_Q(f)} = \frac 1{(q-1)q^{2g}} \sum_{2\alpha+\beta=2g+1} 
\sigma(f;\alpha) \sum_{\deg B=\beta} \leg{B}{f} 
\end{equation*}
where 
\begin{equation*}
 \sigma(f;\alpha):= \sum_{\substack{\deg A=\alpha\\ \gcd(A,f)=1}} \mu(A) \;.
\end{equation*}
\subsection{A sum of M\"obius values} 
Suppose $P$ is a prime of degree $n$, $k\geq 1$ and $\alpha\geq 0$. Set 
$$\sigma_n(\alpha) :=\sigma(P^k;\alpha)=
 \sum_{\substack{ \deg A=\alpha\\ \gcd(A,P^k)=1}} \mu(A) \;.
$$
Since the conditions $\gcd(A,P^k)=1$ and $\gcd(A,P)=1$ are equivalent
for a prime $P$ and any $k\geq 1$, this quantity is independent of $k$;   
the notation anticipates that it depends only on the degree 
$n$ of $P$, as is shown in:   
\begin{lem}\label{lem Mobius} 
i) For $n=1$, 
$$
\sigma_1(0)=1,\quad \sigma_1(\alpha)=1-q \mbox{ for all }
\alpha \geq 1 \;.
$$ 

ii) If $n\geq 2$ then 
\begin{equation*}
\label{mobius sum}
\sigma_n(\alpha) =
\begin{cases} 
1&\alpha =0 \mod n\\ -q&\alpha=1 \mod n \\0& \mbox{ otherwise} 
\end{cases} \;.
\end{equation*}
\end{lem}
\begin{proof}
Since $P$ is prime, 
$$\sigma_n(\alpha) = \sum_{\deg A = \alpha} \mu(A) - 
\sum_{\substack{\deg A = \alpha\\P\mid A}} \mu(A) = 
\sum_{\deg A = \alpha} \mu(A) - \sum_{\deg A_1=\alpha-n} \mu(PA_1) \;.
$$
Now $\mu(PA_1)\neq 0$ only when $A_1$ is coprime to $P$, in which case
$\mu(PA_1)=\mu(P)\mu(A_1)=-\mu(A_1)$. Hence 
$$\sigma_n(\alpha) = \sum_{\deg A = \alpha} \mu(A) + 
\sum_{\substack{\deg A_1 = \alpha-n\\(P,A_1)=1}} \mu(A_1) \;,
$$
that is 
$$ \sigma_n(\alpha)-\sigma_n(\alpha-n) = \sum_{\deg A = \alpha} \mu(A) 
=\begin{cases} 1&\alpha =0\\ -q&\alpha=1 \\0& \alpha \geq 2 \end{cases}
$$
on using 
$$
\sum_{A \mbox{ monic}} \frac{\mu(A)}{|A|^s} = \frac 1{\zeta_q(s)} =1-q^{1-s}
$$
and \eqref{simple zeta}. 
For $n\geq 2$ we get (ii) while for $n=1$ we find that $\sigma_1(0)=1$
and for $\alpha\geq 1$, 
$$\sigma_1(\alpha)=\sigma_1(\alpha-1) = \dots = \sigma_1(1)=-q$$
giving (i). 
\end{proof}

\subsection{The probability that $P\nmid Q$} 

\begin{lem}\label{prob P nmid Q}
 Let $P$ be a prime. Then
$$\ave{\chi_Q(P^2)}= \frac {|P|} {|P|+1} +O(q^{-2 g} ) \;.
$$
\end{lem}
\begin{proof}
Since $P$ is prime, $\chi_Q(P^2)=1$ unless $P$ divides $Q$, that is setting 
$$\iota_P(f):=\begin{cases}1,&P\nmid f\\0,&P\mid f \end{cases} 
$$
we have $\chi_Q(P^2)=\iota_P(Q)$ and thus  by \eqref{ave with mobius}
$$
\ave{\chi_Q(P^2)}= \ave{\iota_P}  = \frac 1{(q-1)q^{2g}} \sum_{\deg A^2 B=2g+1} \mu(A)\iota_P(A^2B) \;.
$$
Since $P$ is prime, $P\nmid A^2B$ if and only if $P\nmid A$ and $P\nmid B$. Hence 
$$\ave{\chi_Q(P^2)} = \frac 1{(q-1)q^{2g}} \sum_{0\leq \alpha \leq g} 
\sum_{\deg A=\alpha, P\nmid A} \mu(A) \sum_{\deg B=2g+1-2\alpha, P\nmid B} 1 \;.
$$
Writing $m:=\deg P$,  
$$ 
\#\{B:\deg B=\beta\quad P\nmid B \} = q^{\beta} 
\cdot \begin{cases} 1,& \mbox{if } m>\beta \\ 1-\frac 1{|P|}, & 
\mbox{if } m\leq \beta \end{cases}
$$
and 
$$ \sum_{\deg A=\alpha, P\nmid A} \mu(A)  = \sigma_m(\alpha)$$
is computed in Lemma~\ref{lem Mobius}. Hence 
\begin{equation*}
\begin{split}
\ave{\chi_Q(P^2)} &= \frac 1{(q-1)q^{2g}} \sum_{0\leq \alpha \leq
  g}\sigma_m(\alpha) q^{2g+1-2\alpha} \cdot 
\begin{cases} 1-\frac 1{|P|},& 0\leq \alpha \leq g-\frac {m-1}2 \\ 1,& g-\frac {m-1}2<\alpha\leq g \end{cases}\\
 &= (1-\frac 1{|P|}) \frac 1{1-\frac 1q} \left( 
\sum_{\alpha =0}^\infty \frac {\sigma_m(\alpha)}{q^{2\alpha}}
+O(q^{-2g}) \right)  \;.
                \end{split}
\end{equation*}
Moreover, inserting the values of $\sigma_m(\alpha)$ given by
Lemma~\ref{lem Mobius} gives  
$$ \sum_{\alpha= 0}^\infty \frac {\sigma_m(\alpha)}{q^{2\alpha}} =\frac{ 1-\frac 1q}{1-\frac 1{|P|^2}}$$
(this is valid both for $m=1$ and $m\geq 2$ !) and hence 
$$
\ave{\chi_Q(P^2)}= (1-\frac 1{|P|}) \frac 1{1-\frac 1q}  \frac{ 1-\frac 1q}{1-\frac 1{|P|^2}} +O(q^{-2g}) 
=\frac {|P|}{|P|+1}+O(q^{-2g}) 
$$
as claimed. 
\end{proof}

\section{Double character sums} 

We consider the double character sum 
\begin{equation*}\label{def of S(beta,n)} 
S(\beta;n):= \sum_{\substack {\deg P=n\\ P \mbox{ prime }}} 
\sum_{\substack{\deg B=\beta\\B \mbox{ monic }}}\leg{B}{P} \;.
\end{equation*}
We may express $S(\beta,n)$ in terms of the coefficients 
$A_P(\beta)=\sum_{\deg B=\beta} \chi_P(B)$ 
of the L-function $\LL(u,\chi_P) = \sum_\beta A_P(\beta) u^\beta$: 
\begin{equation*}
S(\beta;n)  = (-1)^{\frac {q-1}2 \beta n} \sum_{\deg P=n} A_P(\beta) \;,
\end{equation*} 
which follows from the law of quadratic reciprocity \cite{Rosen} : 
If $A$, $B$ are monic then 
$$
\leg{B}{P} = (-1)^{\frac {q-1}2 \deg P \deg B} \leg{P}{B} = 
(-1)^{\frac {q-1}2 \deg P \deg B} \chi_P(B) \;.
$$

Since $A_P(\beta)=0$ for $\beta\geq \deg P$, we find: 
\begin{lem}\label{vanishing of S(beta,n)} 
For $n\leq \beta$ we have  
\begin{equation*}
S(\beta;n) = 0 \;.
\end{equation*}
\end{lem}

\subsection{Duality} 
\begin{prop}
i) If $n$ is odd and $0\leq \beta\leq n-1$ then 
\begin{equation}\label{odd duality} 
S(\beta;n)=q^{\beta-\frac{n-1}2} S( n-1-\beta;n) 
\end{equation}
and 
\begin{equation}\label{S(n-1,n) odd} 
S(n-1;n) =   \pi_q(n) q^{\frac{n-1}2}, \quad n \mbox{ odd} \;.
\end{equation}

ii) If $n$ is even and $1\leq \beta\leq n-2$ then 
\begin{equation}\label{even duality for S} 
S(\beta;n) =  q^{\beta-\frac n2} \left( -S(n-1-\beta;n) + (q-1) \sum_{j=0}^{n-\beta-2}
S(j;n) \right) 
\end{equation}
and 
\begin{equation}\label{S(n-1,n) even} 
S(n-1;n) = -\pi_q(n)q^{\frac{n-2}2}, \quad n \mbox{ even} \;.
\end{equation}
\end{prop}
\begin{proof} 
Assume that $n=\deg P$ is odd. Then
 $\LL(u,\chi_P) = \LL^*(u,\chi_P)$, and so the coefficients
$A_P(\beta) =A_P^*(\beta)$ coincide.   
Therefore the functional equation  in the form \eqref{FE for coeffs} implies  
\begin{equation*}
A_P(\beta) = A_P(n-1-\beta)q^{\beta-\frac{n-1}2} ,\quad n \mbox{ odd}, \quad 0\leq
\beta\leq n-1  \;.
\end{equation*}
Consequently we find that for $n$ odd, 
\begin{equation*}
S(\beta;n)=q^{\beta-\frac{n-1}2} S( n-1-\beta;n) ,\quad n \mbox{ odd}, \quad 0\leq
\beta\leq n-1 \;.
\end{equation*}
In particular we have 
\begin{equation*} 
S(n-1;n) = q^{\frac{n-1}2} S(0,n) = q^{\frac{n-1}2}\pi_q(n), \quad n \mbox{ odd} \;.
\end{equation*}

Next, assume that $n=\deg P$ is even. Then $\LL(u,\chi_P) =
(1-u)\LL^*(u,\chi_P)$, which implies that the coefficients of $\LL(u,\chi_P)$ and $\LL^*(u,\chi_P)$ satisfy 
$$ A_P(\beta)=A_P^*(\beta)-A_P^*(\beta-1) ,\qquad \beta \geq 1 $$ 
and  
\begin{equation}\label{b in terms of A} 
A_P^*(\beta) = A_P(\beta)+A_P(\beta-1)+\dots +A_P(0) \;.
\end{equation}
Moreover 
$$ A_P(0)=A_P^*(0),\quad A_P(n-1) = -A_P^*(n-2) \;.
$$
In particular, since 
\begin{equation*}
A_P^*(0)=1,\quad A_P^*(n-2) = q^{\frac{n-2}2}
\end{equation*}
(see \eqref{FE for coeffs}) we get  
\begin{equation*}\label{AP_(n-1) even} 
A_P(n-1) = -A_P^*(n-2) = - q^{\frac{n-2}2}, \quad n\mbox{ even}
\end{equation*}
so that 
\begin{equation*}
S(n-1;n) = -\pi_q(n)q^{\frac{n-2}2}, \quad n \mbox{ even} \;.
\end{equation*}

The functional equation \eqref{FE for coeffs} implies  
\begin{equation*}
A_P^*(\beta)=A_P^*(n-2-\beta)q^{\beta-\frac{n-2}2} ,\quad 0\leq \beta\leq n-2
\end{equation*}
and hence for $1\leq \beta\leq n-2$
$$ A_P(\beta) =A_P^*(\beta)-A_P^*(\beta-1)
= A_P^*(n-2-\beta)q^{\beta-\frac{n-2}2} - A_P^*(n-1-\beta) q^{\beta-\frac n2} $$
and inserting \eqref{b in terms of A} gives 
\begin{equation*}
A_P(\beta) = q^{\beta-\frac n2} \left( -A_P(n-1-\beta) +(q-1)
\sum_{j=0}^{n-\beta-2} A_P(j) \right) \;.
\end{equation*}
Summing over all primes $P$ of degree $n$ gives 
\begin{equation*}
S(\beta;n) =  q^{\beta-\frac n2} \left( -S(n-1-\beta;n) + (q-1) \sum_{j=0}^{n-\beta-2}
S(j;n) \right)  
\end{equation*}
as claimed. 
\end{proof}

\subsection{An estimate for $S(\beta;n)$}

\begin{lem} \label{estimate for S(beta,n)} 
Suppose $\beta<n$. Then 
\begin{equation}\label{crude estimate S(beta,n)} 
S(\beta;n) =\eta_\beta \pi_q(n) q^{\frac \beta 2} +O(\frac{\beta}{n}
q^{\frac n2+\beta})  
\end{equation}
where $\eta_\beta=1$ for $\beta$ even, and $\eta_\beta=0$ for $\beta $ odd.
\end{lem}
\begin{proof}
We write 
$$
S(\beta;n) = \sum_{\substack{ B=\Box \\ \deg B=\beta}}  \sum_{\deg P=n}
\leg{B}{P} + 
\sum_{\substack{ B \neq \Box \\ \deg B=\beta}}  \sum_{\deg P=n}
\leg{B}{P}
$$
where the squares only occur when $\beta$ is even.

For $B$ not a perfect square, we use the Riemann Hypothesis for curves 
in the form \eqref{estimate on prime sum}:  
$$ \sum_{\deg P=n} \leg{B}{P} \ll   \frac{\deg B}n q^{n/2} \;.
$$
Hence summing over all nonsquare $B$ of degree $\beta$, of which there are
at most $q^\beta$,  gives
$$\sum_{\substack{ B \neq \Box \\ \deg B=\beta}}  \sum_{\deg P=n}
\leg{B}{P} \ll \frac{\beta}{n} q^{\beta+\frac n2} \;.
$$

Assume now that $\beta$ is even. 
For $B=C^2$, we have $P$ and $B$ are coprime since 
$\deg C=\beta/2 <n=\deg P$, and hence $\leg{B}{P}=\leg{C^2}{P}=+1$  
and so the squares, of which there are $q^{\beta/2}$, contribute $\pi_q(n)q^{\beta/2}$. 
This proves \eqref{crude estimate S(beta,n)}. 
\end{proof}

By using duality, \eqref{crude estimate S(beta,n)} 
can be bootstrapped into an improved estimate when $\beta$ is odd: 
\begin{prop}
If $\beta$ is odd and $\beta< n$ then 
\begin{equation}\label{bootstrap}
S(\beta;n)  = -\eta_n \pi_q(n) q^{\beta-\frac n2} + O( q^{n} ) \;.
\end{equation}
\end{prop}
\begin{proof}
Assume $n$ odd with $\beta<n$. 
Then by \eqref{odd duality} for odd $n$,  
$$ S(\beta;n) = q^{\beta-\frac{n-1}2}S(n-1-\beta;n) $$
and inserting the inequality \eqref{crude estimate S(beta,n)} with $\beta$
replaced by $n-1-\beta$ (which is odd in this case) we get 
$$ S(n-1-\beta;n) \ll   q^{\frac n2 +(n-1-\beta)}$$
hence 
$$ S(\beta;n) \ll q^{\beta-\frac{n-1}2}  q^{\frac n2 +(n-1-\beta)} \ll  q^n$$
as claimed. 

Assume $n$ even, with $\beta<n$. Using \eqref{even duality for S} and the
bound \eqref{crude estimate S(beta,n)} gives 
\begin{equation*}
\begin{split}
S(\beta;n) &=  q^{\beta-\frac n2} \left( -S(n-1-\beta;n) + (q-1) \sum_{j=0}^{n-\beta-2}
S(j;n) \right) \\
&= q^{\beta-\frac n2} \left( -\eta_{n-1-\beta} \pi_q(n) q^{\frac{n-1-\beta}2} +(q-1)
\sum_{j=0}^{n-\beta-2}  \eta_j \pi_q(n) q^{\frac j2} \right) \\
&+ O\left( q^{\beta-\frac n2} \sum_{j=0}^{n-1-\beta} \frac jn q^{\frac n2+j}\right)\;.
\end{split}
\end{equation*}
The remainder term is $O(q^n)$. For the main term, we note that
$n-1-\beta=2L$ is even since $\beta$ is odd and $n$ is even, and then we can
write the sum as 
$$
q^{\beta-\frac n2} \pi_q(n) \left(- q^{L}+(q-1)\sum_{l=0}^{L-1} q^l \right)
= -q^{\beta-\frac n2} \pi_q(n) 
$$
which is our claim. 
\end{proof}
 
\section{Proof of Theorem~\ref{one level}} 
The explicit formula  \eqref{explicit formula} says that for $n>0$, 
\begin{equation*}
\tr \Theta_Q^n = -\frac 1{q^{n/2}}\sum_{\deg f=n} \Lambda(f) \chi_Q(f)
\end{equation*} 
the sum over all monic primes powers. We will separately treat the
contributions $\mathcal P_n$ of primes, $\Box_n$ of squares  and $\mathbb H_n$ 
of higher odd powers of primes:
\begin{equation}\label{division of Theta^n} 
\tr \Theta_Q^n = \mathcal P_n + \Box_n + \mathbb H_n \;.  
\end{equation}

\subsection{The contribution of squares } 
When $n$ is even, we have a
contribution to $\tr \Theta_Q^n$ coming from squares of
prime powers (for odd $n$ this term does not appear), which give 
$$
 \Box_n= 
-\frac 1{q^{n/2}} \sum_{\deg h=\frac n2} \Lambda(h) \chi_Q(h^2) \;.
$$
Since $\chi_Q(h^2)=0$ or $1$, we clearly have 
$\Box_n \leq 0$ and 
$$ 
\Box_n \geq -\frac 1{q^{n/2}} \sum_{\deg h=\frac n2} \Lambda(h)  =-1 \;.
$$
by \eqref{equiv simple zeta}. 
Hence the contribution of squares is certainly bounded.

Now for $h=P^k$ a prime power,
\begin{equation}
\ave{\chi_Q(h^2)} = \ave{\chi_Q(P^2)} = 1-\frac{1}{|P|+1} +O(q^{-2g}) \;.
\end{equation} 
by Lemma~\ref{prob P nmid Q}.  
Thus, recalling that $\sum_{\deg h=m} \Lambda(h) = q^m$ 
\eqref{equiv simple zeta}, 
the contribution of squares to the average is 
\begin{equation}
\begin{split}
\ave{\Box_n }  &= 
-1+\frac 1{q^{n/2}} \sum_{\deg P\mid \frac n2}\left( \deg(P)\frac{1}{|P|+1}
+O(q^{-2g}) \right) \\
&= -1+\frac 1{q^{n/2}}  \sum_{\deg P\mid \frac n2}
\frac{\deg(P)}{|P|+1} +O(q^{-2g}) \;.
\end{split}
\end{equation}
In particular, we find that the contribution of squares  to the
average is 
$$  
\ave{\Box_n }= -1+ O (\frac n{q^{n/2}}) +O(q^{-2g}) 
$$
and thus if $n\gg 3\log_q g$ we get
$$
\ave{\Box_n }= -\eta_n(1+ o(\frac 1g)) \;.
$$

\subsection{The contributions of primes} 
The contribution to $\tr \Theta_Q^n$ of primes is 
$$ \mathcal P_n = -\frac n{q^{ n/2}} \sum_{\deg P=n}\chi_Q(P) \;.
$$ 
\begin{prop}
\begin{equation}\label{ave prime} 
\ave{\mathcal P_n} =  -\frac n{(q-1)q^{2g+n/2}} 
\sum_{\substack {\beta+2\alpha = 2g+1\\ \alpha,\beta \geq 0}} \sigma_n( \alpha) S(\beta;n) \;.
  \end{equation}
Moreover, if $n>g$ then 
\begin{equation}\label{when n>g}
\ave{\mathcal P_n} =  -\frac n{(q-1)q^{2g+n/2}} 
\left( S(2g+1;n)- qS(2g-1;n) \right) \;.
\end{equation} 
\end{prop}
\begin{proof}
  Using \eqref{ave with mobius} we have 
\begin{equation*}
\begin{split}
\ave{\mathcal P_n}  
&= -\frac n{(q-1)q^{2g+n/2}} \sum_{\deg P=n}
\sum_{\substack {\beta+2\alpha = 2g+1\\ \alpha,\beta \geq 0}} 
\sigma_n( \alpha) \sum_{\deg B=\beta} \leg{B}{P}    \\
&= -\frac n{(q-1)q^{2g+n/2}} 
\sum_{\substack {\beta+2\alpha = 2g+1\\ \alpha,\beta \geq 0}} \sigma_n( \alpha) S(\beta;n) 
\end{split}
\end{equation*}
which gives \eqref{ave prime}. 

Now assume that $n>g$. Then $\sigma_n(\alpha)\neq 0$ forces
$\alpha=0,1 \mod n$ by Lemma~\ref{lem Mobius}(ii) 
and together with $\alpha \leq g<n$ we must have
$\alpha=0,1$. Hence in \eqref{ave prime} the only nonzero terms are
those with $\alpha=0,1$ which gives \eqref{when n>g}.
\end{proof}

\subsection{Bounding the contribution of primes} 

Assume first that $n\leq g$. In \eqref{ave prime}, if $S(\beta;n)\neq
0$ then  $\beta<n$ by  Lemma~\ref{vanishing of S(beta,n)}. 
For those, we use the   bound
$|S(\beta;n)|\ll \frac{\beta}{n} q^{\beta+n/2}$ of 
Lemma~\ref{estimate for S(beta,n)} and hence  
\begin{equation}\label{primes n<g} 
\ave{\mathcal P_n}   \ll \frac{n}{q^{2g+\frac n2}} 
\sum_{\beta<n} \frac{\beta}{n} q^{n/2+\beta} \ll n q^{n-2g} \leq g q^{-g }
\end{equation}
since $n\leq g$, which vanishes as $g\to \infty$.

For $g<n<2g$, use \eqref{when n>g}, and note that $S(2g\pm 1;n)=0$ by 
Lemma~\ref{vanishing of S(beta,n)}. Hence 
$$\ave{\mathcal P_n}   = 0,\quad g<n<2g\;.
$$

\subsubsection{The case $n=2g$} 
We have $S(2g+1;2g)=0$  by Lemma~\ref{vanishing of S(beta,n)}, and 
$S(2g-1;2g)= -\pi(2g)q^{\frac{2g-2}2}$ by \eqref{S(n-1,n) even}.
Hence 
\begin{equation*}
\begin{split}
\ave{\mathcal P_n } 
&= -\frac {2g}{(q-1)q^{2g+g}} \left( S(2g+1,2g)-qS(2g-1,2g) \right) \\
&= -\frac {2g}{(q-1)q^{2g+g}}q \pi(2g)q^{\frac{2g-2}2} \\
& = -\frac 1{q-1} +O(gq^{-g}) \;.
\end{split} 
\end{equation*}

\subsubsection{The case $2g<n$} 
Here we use \eqref{bootstrap} to find 
\begin{equation*}
\begin{split}
\ave{\mathcal P_n }  &= -\frac {n}{ (q-1)q^{2g+\frac n2}} \left(
S(2g+1;n)-qS(2g-1;n) \right)\\
&=  -\frac {n}{ (q-1)q^{2g+\frac n2}} 
\left(-\eta_n\pi_q(n)q^{2g+1-\frac n2} +q\eta_n\pi_q(n)q^{2g-1-\frac n2} \right)  \\
& +O\left(  \frac n{q^{2g+\frac n2}}  q^{n}\right)\\
&= \eta_n \frac{n\pi_q(n) }{q^n} 
+O\left( n q^{\frac{n}2-2g} \right)  \\
&= \eta_n\left( 1+O(gq^{-g}) \right) + O(n q^{\frac n2-2g}) \;.
\end{split}
\end{equation*}
The main term is asymptotic to $\eta_n$, and 
the remainder is $o(1/g)$  provided 
$$ 2g<n<4g-5\log_q g \;. $$

\subsection{The contribution of higher prime powers} 
The contribution of odd powers of primes $P^d$, $d>1$ odd, $\deg
P^d=n$, is 
$$ 
\mathbb H_n = -\frac 1{q^{\frac n2}} 
\sum_{\substack{d\mid n\\3\leq d \mbox{ odd}}} 
\sum_{\deg P=\frac nd} \frac nd  \chi_Q(P^d) \;.
$$
Since $\chi_Q(P^d)=\chi_Q(P)$ for $d$ odd, the average is 
\begin{equation*}
  \begin{split}
\ave{\mathbb H_n } &= 
-\frac1{(q-1)q^{2g+\frac n2}}\sum_{\substack{d\mid n\\3\leq d \mbox{ odd}}} \frac nd
\sum_{\deg P=\frac nd} \sum_{2\alpha+\beta=2g+1} 
\sigma_{n/d}(\alpha) \sum_{\deg B=\beta} \leg{B}{P}\\
&=    -\frac1{(q-1)q^{2g+\frac n2}}\sum_{\substack{d\mid n\\3\leq d \mbox{
      odd}}} \frac nd 
\sum_{2\alpha+\beta=2g+1} \sigma_{n/d}(\alpha) S(\beta;\frac nd) \;.
  \end{split}
\end{equation*}

In order that $S(\beta;\frac nd)\neq 0$ we need $\beta<n/d$. 
Thus using the bound  
$S(\beta;\frac nd)\ll  q^{\beta+\frac n{2d}}$ 
of \eqref{crude estimate S(beta,n)} 
(recall that $\beta\leq 2g+1$ is odd here) gives 
\begin{equation*}
  \begin{split}
\ave{\mathbb H_n } &\ll 
\frac 1{q^{2g+\frac n2}} 
\sum_{\substack{d\mid n\\3\leq d \mbox{ odd}}} \frac nd 
\sum_{\beta\leq \min(n/d,2g+1)}   q^{\frac{n}{2d}+\beta} \\
& \ll 
 \frac {n} {q^{2g+\frac n2}} 
\sum_{\substack{d\mid n\\3\leq d \mbox{ odd}}}  q^{\frac{ n}{2d}+\min(2g,\frac nd)} \;.
 \end{split}
\end{equation*}
Treating separately the cases $n/3<2g$ and $n/3 \geq 2g$ we see that we have in either case 
\begin{equation}
\ave{\mathbb H_n } \ll g q^{-2g} \;.
\end{equation}


\subsection{Conclusion of the proof} 
We saw that 
$$ 
 \ave{\tr \Theta_Q^n} = \ave{\mathcal P_n } + \ave{\Box_n } +
 \ave{\mathbb H_n }  
$$
with the individual terms giving 
$$
\ave{\mathcal P_n } = 
\begin{cases} 
 O(g q^{-g}), & 0<n<2g\\ -\frac 1{q-1}+O(gq^{-g}), & n=2g \\ 
\eta_n+ O(n q^{n/2-2g}), & 2g<n
\end{cases} \;,
$$
$$ 
\quad \ave{\Box_n } =-\eta_n + \eta_n\frac 1{q^{n/2}} 
\sum_{\deg P\mid \frac n2} \frac{\deg P}{|P|+1} +O(q^{-2g}) \;,
$$ 
and   
$$\ave{\mathbb H_n } =O(g q^{-2g}) \;.$$
Putting these together gives Theorem~\ref{one level}. In particular 
$$ \ave{\tr \Theta_Q^n} = 
\left\{ 
\begin{matrix} 
-\eta_n,& 3\log_q g<n<2g \\  \\ -1-\frac 1{q-1},& n=2g \\  \\ 
0,& 2g<n<4g-8\log_q g
\end{matrix}\right\}   +o(\frac 1g)  \;.
 $$


\section{The product of two traces} \label{sec:products} 
 Using the methods of this paper, on can also compute mean values of products of traces. 
For the product of two traces, the results can be stated as follows:

Assume $\min(m,n)\gg \log g$ and $m+n \leq 4g-100 \log_q g$. 
Then  

i) If $m=n$ then 
\begin{equation}
  \ave{|\tr \Theta_Q^n|^2}  \sim  
  \begin{cases}
    n+\eta_n, & n<g\\ n+\eta_n +\frac 1{q-1},& n=g\\
n-1+\eta_n, & g<n<2g-50\log_q g
  \end{cases}
\end{equation}

ii) If $m<n$ then for ``generic'' values of $(m,n)$ we have
\begin{equation}\label{generic m,n}
  \ave{\tr \Theta_Q^m \tr \Theta_Q^n} \sim  
  \begin{cases}
    \eta_m \eta_n,& m+n<2g \\ 
\eta_m\eta_n-\eta_{m+n}, &  n<2g, \quad m+n>2g \\
-\eta_{m+n},& n>2g, \quad n-m<2g\\
0,& n-m>2g\\
  \end{cases}
\end{equation}
while on ``exceptional'' lines we have 
\begin{equation}\label{exceptional lines} 
  \ave{\tr \Theta_Q^m \tr \Theta_Q^n} \sim  
  \begin{cases}
\eta_m\eta_n+\frac 1{q-1}, & m+n=2g \\
\eta_m\eta_n- \eta_{m+n}+ \eta_m\frac 1{q-1},& n=2g \\
-\frac q{q-1} \eta_{m+n}, & n-m=2g
\end{cases}
\end{equation}

The expected values for the symplectic group are \cite{DS,DE,HR,KO}: 

i) If $m=n$ then 
\begin{equation}\label{RMT diagonal pairs}
  \int_{\USp(2g)}|\tr U^n|^2 dU  = 
  \begin{cases}
    n+\eta_n,& 1\leq n \leq g \\ n-1+\eta_n,& g+1\leq n \leq 2g \\ 2g,& n>2g
  \end{cases}
\end{equation}

ii) If $1\leq m<n$
\begin{equation}\label{RMT generic pairs}
  \int_{\USp(2g)}\tr U^m\tr U^n dU  = 
  \begin{cases}
    \eta_m\eta_n,& m+n \leq 2g \\ 
\eta_m\eta_n-\eta_{m+n},& m<n\leq 2g, \quad m+n>2g\\
-\eta_{m+n},& n>2g,\quad n-m\leq 2g \\
0,& n-m>2g
  \end{cases}
\end{equation}

 Comparing \eqref{RMT diagonal pairs}, \eqref{RMT generic pairs},  and \eqref{generic m,n}, \eqref{exceptional lines} 
we find that if $m=\min(m,n)\gg \log_q g$ and $m+n<4g-100  \log_q g$ then for ``generic'' values of $(m,n)$, that
 is if $n,m\neq 2g$ and  $|n \pm m| \neq 2g$ then 
 \begin{equation}
\ave{\tr \Theta_Q^m \tr \Theta_Q^n} \sim   \int_{\USp(2g)}\tr U^m\tr U^n dU
 \end{equation}
while on the lines $n,m=2g$, $|n\pm m|=2g$ the difference between the averages over $\mathcal H_{2g+1}$ 
and $\USp(2g)$ is bounded by 
\begin{equation*}
|  \ave{\tr \Theta_Q^m \tr \Theta_Q^n} - \int_{\USp(2g)}\tr U^m\tr U^n dU |
\leq \frac1{q-1} +o(1), \qquad g\to \infty \;.
\end{equation*}

These results can be use to study the two-level density, compare \cite{Rubinstein, Gao}.

\end{document}